\documentclass[leqno,10pt]{amsart}

\usepackage{amsmath,amsfonts,amsthm}
\usepackage{amssymb}
\usepackage[all]{xy}
\usepackage{enumerate}
\usepackage{mathrsfs}	
\usepackage{pstricks}
\usepackage{pst-3dplot,pstricks-add}
\usepackage{graphicx}

\theoremstyle{plain}
\newtheorem{thm}{Theorem}[section]

\theoremstyle{definition}
\newtheorem{defn}[thm]{Definition}
\newtheorem{exmp}[thm]{Example}

\def\NN{\mathbb{N}}

\def\kk{\mathbb{K}}
\def\PP{\mathbb{P}}

\def\Sc{\mathcal{S}}

\def\Zc{\mathcal{Z}}

\def\Nc{\mathcal{N}}

\def\Tc{\mathscr{T}}

\def\conv{\mathrm{conv}}

\def\im{\mathrm{im}}

\def\dim{\mathrm{dim}}

\def\dto{\mathfrak{\dashrightarrow}}

\title[A package for computing implicit equations from toric surfaces]
{A package for computing implicit equations of parametrizations from toric surfaces}

\author{Nicol\'{a}s Botbol}
\address{Departamento de Matem\'atica\\
FCEN, Universidad de Buenos Aires, Argentina \\
\& Institut de Math\'ematiques de Jussieu \\
Universit\'e de P. et M. Curie, Paris VI, France \\
E-mail address: nbotbol@dm.uba.ar
}

\author{Marc Dohm}

\thanks{The authors were partially supported by the project ECOS-Sud A06E04. NB was 
partially supported by UBACYT X064, CONICET PIP 112-200801-00483 and ANPCyT PICT
20569, Argentina. MD was partially supported by the project GALAAD, INRIA Sophia Antipolis, France.}

\begin{document}

\begin{abstract}
In this paper we present an algorithm for computing a \textit{matrix representation} for a surface in $\PP^3$ parametrized over a 2-dimensional toric variety $\Tc$. This algorithm follows the ideas of \cite{BDD08} and it was implemented in Macaulay2 \cite{M2}. We showed in \cite{BDD08} that such a surface can be represented by a matrix of linear syzygies if the base points are finite in number and form locally a complete intersection, and in \cite{Bot09} we generalized this to the case where the base locus is not necessarily a local complete intersection. The key point consists in exploiting the sparse structure of the parametrization, which allows us to obtain significantly smaller matrices than in the homogeneous case.
\end{abstract}


\maketitle


\section{Introduction}

Let $\Tc$ be a two-dimensional projective toric variety, and let $f:\Tc \dto \PP^3$ be a generically finite rational map. Hence, $\Sc:=\overline{\im(f)}\subset \PP^3$ is a hypersurface. In \cite{BDD08} and \cite{Bot09} we showed how to compute an implicit equation for $\Sc$, assuming that the base locus $X$ of $f$ is finite and locally an almost complete intersection. The work in \cite{BDD08} and \cite{Bot09} is a further generalization of the results in \cite{BuJo03, BC05, Ch06, BD07} on implicitization of rational hypersurfaces via approximation complexes. 

We showed in \cite{BDD08} how to compute a symbolic matrix of linear syzygies $M$, called \textit{representation matrix} of $\Sc$, with the property that, given a point $p\in \PP^3$, the rank of $M(p)$ drops if $p$ lies in the surface $\Sc$. When the base locus $X$ is locally a complete intersection, we get that the rank of $M(p)$ drops if and only if $p$ lies in the surface $S$. 

We begin by recalling the notion of a representation matrix.

\begin{defn}
Let $\Sc \subset \PP^n$ be a hypersurface. A matrix $M$ with entries in the polynomial ring $\kk[T_0,\ldots,T_n]$ is called a \textit{representation matrix} of $\Sc$ if it is generically of full rank and if the rank of $M$ evaluated in a point $p$ of $\PP^n$ drops if and only if the point $p$ lies on $\Sc$.
\end{defn}

It follows immediately that a matrix $M$ represents $\Sc$ if and only if the greatest common divisor $D$ of all its minors of maximal size is a power of a homogeneous implicit equation $F \in \kk[T_0,\ldots,T_n]$ of $\Sc$. When the base locus is locally an almost complete intersection, we can construct a a matrix $M$ such that $D$ factors as $D=F^\delta G$ where $\delta\in \NN$ and $G \in \kk[T_0,\ldots,T_n]$. In \cite[Sec. 3.2]{Bot09}, we gave a description of the surface $(D=0)$ 
In this paper we present an implementation of our results in the computer aided software Macaulay2 \cite{M2}.From a practical point of view our results are a major improvement, as it makes the method applicable for a wider range of parametrizations (for example, by avoiding unnecessary base points with bad properties) and
leads to significantly smaller representation matrices.

There are several advantages of this perspective. The method works in a very general setting and makes only minimal assumptions on the parametrization. In particular, as we have mentioned, it works well in the presence of ``nice'' base points. Unlike the method of toric resultants (cf. for example \cite{KD06}), we do not have to extract a maximal minor of unknown size, since the matrices are generically of full rank. The monomial structure of the parametrization is exploited, in \cite{Bot09} we defined
\begin{defn}\label{defNf}
Given a list of polynomials $f_0,\hdots,f_r$, we define 
\[
 \Nc(f_0,\hdots,f_r):=\conv(\bigcup_{i=0}^r \Nc(f_i)),
\]
the convex hull of the union of the Newton polytopes of $f_i$, and we will refer to this polytope as the \textit{Newton polytope} of the list $f_0,\hdots,f_r$. When $f$ denotes the rational map defining $\Sc$, we will write $\Nc(f):=\Nc(f_1,f_2,f_3,f_4)$, and we will refer to it as the Newton polytope of $f$. 
\end{defn}
In this terms, in our algorithm we fully exploit the structure of $\Nc(f)$, so one obtains much better results for sparse parametrizations, both in terms of computation time and in terms of the size of the representation matrix. Moreover, it subsumes the known method of approximation complexes in the case of dense homogeneous parametrizations. One important point is that representation matrices can be efficiently constructed by solving a linear system of relatively small size (in our case $\dim_\kk(A_{\nu+d})$ equations in $4 \dim_\kk(A_\nu)$ variables). This means that their computation is much faster than the computation of the implicit equation and they are thus an interesting alternative as an implicit representation of the surface. 

On the other hand, there are a few disadvantages. Unlike with the toric resultant or the method of moving planes and surfaces, the matrix representations are not square and the matrices involved are generally bigger than with the me\-thod of moving planes and surfaces. It is important to remark that those disadvantages are inherent to the choice of the method: A square matrix built from linear syzygies does not exist in general and it is an automatic consequence that if one only uses linear syzygies to construct the matrix, it has to be bigger than a matrix which also uses entries of higher degree (see \cite{BCS09}). The choice of the method to use depends very much on the given parametrization and on what one needs to do with the matrix representation.

\section{Example}\label{sec:examples}

\begin{exmp} \label{interestingexample}
Here we give an example, where we fully exploit the structure of $\Nc(f)$. Take $(f_1,f_2,f_3,f_4) = (st^6+2,st^5-3st^3,st^4+5s^2t^6,2+s^2t^6)$. This is a very sparse parametrization, and we have in this case, there is no smaller lattice homothety of $\Nc(f)$ (cf. \cite{BDD08, Bot09} for a wider discussion on this subject). The coordinate ring is $A=\kk[X_0,\ldots,X_5]/J$, where $J=(X_3^2-X_2X_4, X_2X_3-X_1X_4, X_2^2-X_1X_3, X_1^2-X_0X_5)$ and the new base-point-free parametrization $g$ is given by $(g_1,g_2,g_3,g_4)=(2X_0+X_4,-3X_1+X_3, X_2+5X_5, 2X_0+X_5)$. The Newton polytope looks as follows.
\begin{center}
  \includegraphics[scale=0.67]{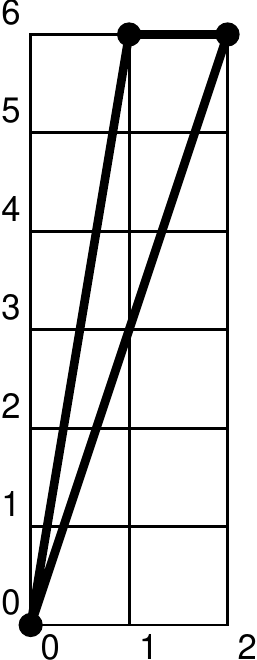}
\end{center}

For $\nu_0=2d=2$ we can compute the matrix of the first map of the graded piece of degree $\nu_0$ of the approximation complex of cycles $(\Zc_\bullet)_{\nu_0}$ (see for example \cite[Sec 3.1]{BDD08}), which is a $17 \times 34$-matrix. 
The greatest common divisor of the $17$-minors of this matrix is the homogeneous implicit equation of the surface; it is of degree 6 in the variables 
{\small
\[
\begin{array}{rl}
 T_1,\ldots,T_4: & 2809T_1^2T_2^4 + 124002T_2^6 - 5618T_1^3T_2^2T_3 + 66816T_1T_2^4T_3 + 2809T_1^4T_3^2\vspace{0.12cm} \\ \vspace{0.12cm}
 &- 50580T_1^2T_2^2T_3^2  + 86976T_2^4 T_3^2 + 212T_1^3T_3^3  - 14210T_1T_2^2T_3^3  + 3078T_1^2 T_3^4 \\  \vspace{0.12cm}
 & + 13632T_2^2 T_3^4  + 116T_1T_3^5 + 841T_3^6  + 14045T_1^3 T_2^2 T_4 - 169849T_1T_2^4 T_4 \\ \vspace{0.12cm}
 & -14045T_1^4 T_3T_4 + 261327T_1^2 T_2^2 T_3T_4 - 468288T_2^4 T_3T_4 - 7208T_1^3 T_3^2 T_4 \\ \vspace{0.12cm}
 & + 157155T_1T_2^2 T_3^3 T_4 - 31098T_1^2 T_3^3 T_4 - 129215T_2^2 T_3^3 T_4 - 4528T_1T_3^4 T_4  \\ \vspace{0.12cm}
 & - 12673T_3^5 T_4 - 16695T_1^2 T_2^2 T_4^2  + 169600T_2^4 T_4^2  + 30740T_1^3 T_3T_4^2 \\ \vspace{0.12cm}
 & - 433384T_1T_2^2 T_3T_4^2 + 82434T_1^2 T_3^2 T_4^2  + 269745T_2^2 T_3^2 T_4^2  + 36696T_1T_3^3 T_4^2 \\ \vspace{0.12cm}
 &  + 63946T_3^4 T_4^2  + 2775T_1T_2^2 T_4^3  - 19470T_1^2 T_3T_4^4  + 177675T_2^2 T_3T_4^3   \\  \vspace{0.12cm}
 &- 85360T_1T_3^2 T_4^3  - 109490T_3^3 T_4^3  - 125T_2^2 T_4^4  + 2900T_1T_3T_4^4  + 7325T_3^2 T_4^4    \\ 
 &- 125T_3T_4^5
\end{array} 
\]
}

In this example we could have considered the parametrization as a bihomogeneous map either of bidegree $(2,6)$ or of
bidegree $(1,3)$, i.e. we could have chosen the corresponding rectangles instead
of $\Nc(f)$. This leads to a more complicated coordinate ring in $20$ (resp.~$7$) variables and $160$ (resp.~$15$) generators of $J$ and to bigger matrices
(of size $21 \times 34$ in both cases). Even more importantly, the parametrizations will have a non-LCI base point and the matrices do not represent the implicit equation but a multiple of it (of degree $9$). Instead, if we consider the map as a homogeneous map of degree $8$, the results are even worse: For $\nu_0 = 6$, the $28 \times 35$-matrix $M_{\nu_0}$ represents a multiple of the implicit equation of degree $21$.

To sum up, in this example the toric version of the method of approximation complexes works well, whereas it fails over $\PP^1 \times \PP^1$ and $\PP^2$. This shows that the extension of the method to toric varieties really is a generalization and makes the method applicable to a larger class of parametrizations.\medskip

Interestingly, we can even do better than with $\Nc(f)$ by choosing a smaller polytope. The philosophy is that the choice of the optimal polytope is a compromise between two criteria: keep the simplicity of the polytope in order not to make the the ring $A$ too complicated, and respect the sparseness of the parametrization (i.e. keep the polytope close to the Newton polytope) so that no base points appear which are not local complete intersections.

So let us repeat the same example with another polytope $Q$, which is small enough to reduce the size of the matrix but which only adds well-behaved (i.e. local complete intersection) base points:
\begin{center}
\includegraphics[scale=0.7]{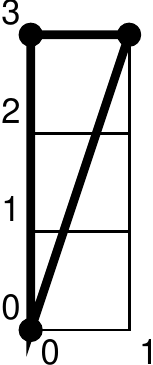}
\end{center}
 The Newton polytope $\Nc(f)$ is contained in $2 \cdot Q$, so the parametrization
will factor through the toric variety associated to $Q$, more precisely we obtain
a new parametrization defined by
 $$(g_1,g_2,g_3,g_4)=(2X_0^2+X_3X_4,-3X_0X_4+X_2X_4, X_1X_4+5X_4^2,2X_0^2+X_4^2)$$
over the coordinate ring $A=\kk[X_0,\ldots,X_4]/J$ with $J=(X_2^2-X_1X_3, X_1X_2-X_0X_3, X_1^2-X_0X_2)$. The optimal bound is $\nu_0=2$ and in this degree the implicit equation is represented directly without extraneous factors by a $12 \times 19$-matrix, which is smaller than the $17 \times 34$ we had before. 
\end{exmp}

\medskip


\section{Implementation in Macaulay2}
In this section we show how to compute a matrix representation and the implicit equation
with the method developed in \cite{BDD08} and \cite{Bot09}, using the computer algebra system Macaulay2 \cite{M2}.
As it is probably the most interesting case from a practical point of view, we restrict our computations to parametrizations of a toric surface.
However, the method can be adapted to the $n$-dimensional toric case. 
Moreover, we are not claiming that our implementation is optimized for efficiency; anyone
trying to implement the method to solve computationally involved examples is well-advised to
give more ample consideration to this issue. For example, in the toric case there are better suited
software systems to compute the generators of the toric ideal $J$, see \cite{4ti2}.

First we load the package ``Maximal minors\footnote{The package ``maxminor.m2" for Macaulay2 can be downloaded from the webpage \texttt{http://mate.dm.uba.ar/\textasciitilde nbotbol/maxminor.m2}.}''
{\small\begin{verbatim}
i1 : load "maxminor.m2"
\end{verbatim}}

Let us start by defining the parametrization $f$ given by $(f_1,\ldots,f_4)$.

{\small\begin{verbatim}
i2 : S=QQ[s,u,t,v];
i3 : e1=2;
i4 : e2=6;
i5 : f1=s*u*t^6+2*u^2*v^6
          6     2 6
o5 = s*u*t  + 2u v
i6 : f2=s*u*t^5*v-3*s*u*t^3*v^3
          5          3 3
o6 = s*u*t v - 3s*u*t v
i7 : f3=s*u*t^4*v^2+5*s^2*t^6
       2 6        4 2
o7 = 5s t  + s*u*t v
i8 : f4=2*u^2*v^6+s^2*t^6
      2 6     2 6
o8 = s t  + 2u v

\end{verbatim}}

We construct the matrix asasociated to the polynomials and we relabel them in order to be able to automatize some procedures.
{\small\begin{verbatim}

i9 : F=matrix{{f1,f2,f3,f4}}
o9 = | sut6+2u2v6 sut5v-3sut3v3 5s2t6+sut4v2 s2t6+2u2v6 |
             1       4
o9 : Matrix S  <--- S
i10 : f_1=f1;
i11 : f_2=f2;
i12 : f_3=f3;
i13 : f_4=f4;
\end{verbatim}}

We define the associated affine polynomials \texttt{FF$\_$i} by specializing the variables $u$ and $v$ to $1$.
{\small\begin{verbatim}
i14 : for i from 1 to 4 do (
        FF_i=substitute(f_i,{u=>1,v=>1});
      )
\end{verbatim}}

We just change the polynomials \texttt{FF$\_$i} to the new ring $S2$.

{\small\begin{verbatim}
i15 : S2=QQ[s,t]
o15 = S2
o15 : PolynomialRing
i16 : for i from 1 to 4 do (
        FF_i=sub(FF_i,S2);
      )
\end{verbatim}}

The reader can experiment with the implementation simply by changing the definition of the polynomials
and their degrees, the rest of the code being identical. We first set up the list $st$ of monomials $s^it^j$ of bidegree $(e'_1,e'_2)$. In the toric case,
this list should only contain the monomials corresponding to points in the Newton polytope $\Nc'(f)$.
{\small\begin{verbatim}
i17 : use S;
i18 : st={};
i19 : for i from 1 to 4 do (
        st=join(st,flatten entries monomials f_i); 
      )
i20 : l=length(st)-1;
i21 : k=gcd(e1,e2)
o21 = 2
\end{verbatim}}

 
We compute the ideal $J$ and the quotient ring $A$. This is done
by a Gr\"obner basis computation which works well for examples of
small degree, but which should be replaced by a matrix formula
in more complicated examples. In the toric case,
there exist specialized software systems such as \cite{4ti2} to compute the ideal $J$.
{\small\begin{verbatim}
i24 : SX=QQ[s,u,t,v,w,x_0..x_l,MonomialOrder=>Eliminate 5]
o24 = SX
o24 : PolynomialRing
i25 : X={};
i26 : st=matrix {st};
              1       8
o26 : Matrix S  <--- S
i27 : F=sub(F,SX)
o27 = | sut6+2u2v6 sut5v-3sut3v3 5s2t6+sut4v2 s2t6+2u2v6 |
               1        4
o27 : Matrix SX  <--- SX
i28 : st=sub(st,SX)
o28 = | sut6 u2v6 sut5v sut3v3 s2t6 sut4v2 s2t6 u2v6 |
               1        8
o28 : Matrix SX  <--- SX
i29 : te=1;
i30 : for i from 0 to l do ( te=te*x_i )
i31 : J=ideal(1-w*te)
o31 = ideal(- w*x x x x x x x x  + 1)
                 0 1 2 3 4 5 6 7
o31 : Ideal of SX
i32 : for i from 0 to l do (
          J=J+ideal (x_i - st_(0,i))
          )
i33 : J= selectInSubring(1,gens gb J)
o33 = | x_4-x_6 x_1-x_7 x_3^2-x_6x_7 x_2x_3-x_5^2 x_0x_3-x_2x_5 
      ---------------------------------------------------------
      x_2^2-x_0x_5 x_5^3-x_0x_6x_7 x_3x_5^2-x_2x_6x_7 |
               1        8
o33 : Matrix SX  <--- SX
i34 : R=QQ[x_0..x_l]
o34 = R
o34 : PolynomialRing
i35 : J=sub(J,R)
o35 = | x_4-x_6 x_1-x_7 x_3^2-x_6x_7 x_2x_3-x_5^2 x_0x_3-x_2x_5 
      ---------------------------------------------------------
      x_2^2-x_0x_5 x_5^3-x_0x_6x_7 x_3x_5^2-x_2x_6x_7 |
              1       8
o35 : Matrix R  <--- R
i36 : A=R/ideal(J)  
o36 = A
o36 : QuotientRing
\end{verbatim}}

Next, we set up the list $ST$ of monomials $s^it^j$ of bidegree $(e_1,e_2)$ and the list $X$ of
the corresponding elements of the quotient ring $A$. In the toric case,
this list should only contain the monomials corresponding to points in the Newton polytope $\Nc(f)$.
{\small\begin{verbatim}
i37 : use SX
o37 = SX
o37 : PolynomialRing
i38 :   ST={};
i39 :   X={};
i40 :   for i from 0 to l do (
          ST=append(ST,st_(0,i)); 
          X=append(X,x_i);
        )
\end{verbatim}}
We can now define the new parametrization $g$ by the polynomials $g_1,\ldots,g_4$. 
 {\small\begin{verbatim}
i41 : X=matrix {X};
               1        8
o41 : Matrix SX  <--- SX
i42 : X=sub(X,SX)
o42 = | x_0 x_1 x_2 x_3 x_4 x_5 x_6 x_7 |
               1        8
o42 : Matrix SX  <--- SXX=matrix {X};
i43 : (M,C)=coefficients(F,Variables=>{s_SX,u_SX,t_SX,v_SX},Monomials=>ST)
o43 = (| sut6 u2v6 sut5v sut3v3 s2t6 sut4v2 s2t6 u2v6 |, {8} | 1 0  0 0 |)
                                                         {8} | 0 0  0 0 |
                                                         {8} | 0 1  0 0 |
                                                         {8} | 0 -3 0 0 |
                                                         {8} | 0 0  0 0 |
                                                         {8} | 0 0  1 0 |
                                                         {8} | 0 0  5 1 |
                                                         {8} | 2 0  0 2 |
o43 : Sequence
i44 : G=X*C
o44 = | x_0+2x_7 x_2-3x_3 x_5+5x_6 x_6+2x_7 |
               1        4
o44 : Matrix SX  <--- SX
i45 : G=matrix{{G_(0,0),G_(0,1),G_(0,2),G_(0,3)}}
o45 = | x_0+2x_7 x_2-3x_3 x_5+5x_6 x_6+2x_7 |
               1        4
o45 : Matrix SX  <--- SX
i46 : G=sub(G,A)
o46 = | x_0+2x_7 x_2-3x_3 x_5+5x_6 x_6+2x_7 |
              1       4
o46 : Matrix A  <--- A
\end{verbatim}}
In the following, we construct the matrix representation $M$. For simplicity, we
compute the whole module $\Zc_1$, which is not necessary as we only need the graded
part $(\Zc_1)_{\nu_0}$. In complicated examples, one should compute only this graded part
by directly solving a linear system in degree $\nu_0$. Remark that the best
bound $\mathrm{nu}= \nu_0$ depends on the parametrization.
{\small\begin{verbatim}
i47 : use A
o47 = A
o47 : QuotientRing
i48 : Z0=A^1;
i49 : Z1=kernel koszul(1,G);
i50 : Z2=kernel koszul(2,G);
i51 : Z3=kernel koszul(3,G);
i52 : nu=-1
o52 = -1
i53 : d=1
o53 = 1
i54 : hfnu = 1
o54 = 1
i55 : while hfnu != 0 do (
      nu=nu+1;
      hfZ0nu = hilbertFunction(nu,Z0);
      hfZ1nu = hilbertFunction(nu+d,Z1);
      hfZ2nu = hilbertFunction(nu+2*d,Z2);
      hfZ3nu = hilbertFunction(nu+3*d,Z3);
      hfnu = hfZ0nu-hfZ1nu+hfZ2nu-hfZ3nu;
       );
i56 : nu
o56 = 2
i57 : hfZ0nu
o57 = 17
i58 : hfZ1nu
o58 = 34
i59 : hfZ2nu
o59 = 23
i60 : hfZ3nu
o60 = 6
i61 : hfnu
o61 = 0

i62 : hilbertFunction(nu+d,Z1)-2*hilbertFunction(nu+2*d,Z2)+
      3*hilbertFunction(nu+3*d,Z3)
o62 = 6
i63 : GG=ideal G
o63 = ideal (x  + 2x , x  - 3x , x  + 5x , x  + 2x )
              0     7   2     3   5     6   6     7
o63 : Ideal of A
i64 : GGsat=saturate(GG, ideal (x_0..x_l))
o64 = ideal 1
o64 : Ideal of A
i65 : degrees gens GGsat
o65 = {{{0}}, {{0}}}
o65 : List
i66 : H=GGsat/GG
o66 = subquotient (| 1 |, | x_0+2x_7 x_2-3x_3 x_5+5x_6 x_6+2x_7 |)
                                1
o66 : A-module, subquotient of A
i67 : degrees gens H
o67 = {{{0}}, {{0}}}
o67 : List
i68 : S=A[T1,T2,T3,T4]
o68 = S
o68 : PolynomialRing
i69 : G=sub(G,S);
              1       4
o69 : Matrix S  <--- S
i70 : Z1nu=super basis(nu+d,Z1); 
              4       34
o70 : Matrix A  <--- A
i71 : Tnu=matrix{{T1,T2,T3,T4}}*substitute(Z1nu,S);
              1       34
o71 : Matrix S  <--- S
i72 : 
      lll=matrix {{x_0..x_l}}
o72 = | x_0 x_7 x_2 x_3 x_6 x_5 x_6 x_7 |
              1       8
o72 : Matrix A  <--- A
i73 : lll=sub(lll,S)
o73 = | x_0 x_7 x_2 x_3 x_6 x_5 x_6 x_7 |
              1       8
o73 : Matrix S  <--- S
i74 : ll={}
o74 = {}
o74 : List
i75 : for i from 0 to l do { ll=append(ll,lll_(0,i)) }
i76 : (m,M)=coefficients(Tnu,Variables=>ll,Monomials=>substitute(basis(nu,A),S));
i77 : M;  
              17       34
o77 : Matrix S   <--- S
\end{verbatim}}
The matrix $M$ is the desired matrix representation of the surface $\Sc$.

We can continue by computing the implicit equation and verifying the result by substituting
{\small\begin{verbatim}
 
i78 : T=QQ[T1,T2,T3,T4]
o78 = T
o78 : PolynomialRing
i79 : ListofTand0 ={T1,T2,T3,T4}
o79 = {T1, T2, T3, T4}
o79 : List
i80 : for i from 0 to l do { ListofTand0=append(ListofTand0,0) };
i81 : p=map(T,S,ListofTand0) 
o81 = map(T,S,{T1, T2, T3, T4, 0, 0, 0, 0, 0, 0, 0, 0})
o81 : RingMap T <--- S
i82 : N=MaxCol(p(M)); 
              17       17
o82 : Matrix T   <--- T
i83 : Eq=det(N); factor Eq
\end{verbatim}}
We verify the result by substituting on the computed equation, the polynomials $f_1$ to $f_4$.

{\small\begin{verbatim}
i85 :use S; Eq=sub(Eq,S)
o86 : S
i87 : sub(Eq,{T1=>G_(0,0),T2=>G_(0,1),T3=>G_(0,2),T4=>G_(0,3)})
o87 = 0
\end{verbatim}}


{\small
\def\cprime{$'$}
}


\end{document}